\documentclass{amsart}
\usepackage{setspace, amsmath, amsthm, amssymb, amsfonts, amscd, epic, graphicx, ulem, dsfont}
\usepackage[T1]{fontenc}
\usepackage{multirow}
\usepackage{bbm}

\newtheorem{lem}{Lemma}
\newtheorem{thm}{Theorem}
\newtheorem*{un-thm}{Theorem}

\newtheorem{pro}{Proposition}
\newtheorem{cor}{Corollary}

\theoremstyle{abstract}

\theoremstyle{remark}
\newtheorem*{rema}{Remark}

\newtheorem*{exa}{\textbf{Example 1}}

\theoremstyle{definition}
\newtheorem*{defn}{Definition}

\newcommand{\R}{\mathbb{R}}

\makeatletter \@namedef{subjclassname@2010}{
  \textup{2010} Mathematics Subject Classification}
\makeatother

\begin{document}

\title[Unbounded Products of Operators]{Unbounded Products of Operators and Connections to Dirac-Type Operators}
\author{Karl Gustafson$^1$ and
Mohammed Hichem Mortad$^2$*}
\address{$^1$Department of Mathematics,
University of Colorado, Boulder, CO, USA.}
\email{gustafs@euclid.colorado.edu} \address{$^2$Department of
Mathematics, University of Oran, B.P. 1524, El Menouar, Oran 31000.
Algeria.\newline {\bf Mailing address (for the second author)}:
\newline Dr Mohammed Hichem Mortad \newline BP 7085
Es-Seddikia\newline Oran
\newline 31013 \newline Algeria}

\email{mhmortad@gmail.com, mortad@univ-oran.dz.}

\begin{abstract}
Let $A$ and $B$ be two densely defined unbounded closeable operators
in a Hilbert space such that their unbounded operator products $AB$
and $BA$ are also densely defined. Then all four operators possess
adjoints and we obtain new inclusion bounds for the operator product
closures $\overline{A}~\overline{B}$ and $\overline{AB}$ in terms of
new relations among the operator adjoints. These in turn lead to
sharpened understandings for when products of unbounded self-adjoint
and unbounded normal operators are self-adjoint and normal. They
also clarify certain operator-product issues for Dirac operators.
\end{abstract}



\thanks{* Corresponding author}
\maketitle

\section{Introduction}Let $A$ and $B$ be two densely
defined unbounded closeable operators in a Hilbert space such that
their unbounded operator product domain $D(AB)$, and often, $D(BA)$,
are also dense.  These operators then possess adjoints and we obtain
new results for relations among them.  Our goal is to seek somewhat
finer detail than we have seen in the established literature.  We
shall refrain inasmuch as is possible from any supplementary
assumptions of boundedness, positivity and closed operator ranges.
In particular, this allows our results to be applicable over the
continuous spectrum. See Remark 4.1.

Section 2 investigates the adjoint product $(AB)^*$ and then
establishes inclusion bounds for the operator closures
$\overline{A}~\overline{B}$ and $\overline{AB}$.  Section 3
specializes to $A$ and $B$ self-adjoint and concludes with $A$ and
$B$ normal. Although we are fully cognizant of the power and
naturalness of an assumption of commutativity in such
considerations, a theme there is to also examine its necessity.
Section 4 elaborates key issues and gives some related literature
from which further relevant bibliography may be obtained. In
particular, we clarify and improve certain important results used in
investigations of Dirac operators, see \cite{Gesztesy,Schmincke}.
Thereby, this paper also provides a brief up-to-date account of
these important problems.  An interested reader may therefore wish
to consult Section 4 while reading the earlier sections of this
paper.

\section{Unbounded Product Adjoints}

Recall that a generally unbounded operator $T$ in a Hilbert space
possesses a unique adjoint $T^*$ iff $D(T)$ is dense.  Then $T$ is
closeable with minimal closed extension  $\overline{T}$ iff $D(T^*)$
is dense, and then $(\overline{T})^* = T^*$.  When $A$ and $B$ and
$AB$ are densely defined, then $(AB)^* \supset B^* A^*$, whether or
not the latter adjoint operators are densely defined.  One has $T =
T^{**}$ for closed densely defined operators $T$.  Moreover, if $T$
is just densely defined and closeable, still one has $\overline{T} =
T^{**}$.  Recall adjoint operators are always closed.

In the following we will always assume that our starting $A$ and $B$ are
densely defined so that they have adjoints.

A central question is when $(AB)^* = B^* A^*$.  Such is the case for
everywhere defined bounded operators $A$ and $B$ in $B(H)$. One
needs only $A \in B(H)$ for that equality.  Another central question
is when $AB$ is closed.  When $A$ is closed and $B$ is in $B(H)$,
then $AB$ is closed.  Also when $A$ is invertible with $A^{-1}$ in
$B(H)$, $AB$ is closed.  Some further known results may be found in
the books \cite{Gold-Book-1966,Kato-Book,Riesz-Nagy,Weidmann}. We
digress a bit to say that similar works exist for sums of linear
operators, see e.g.
\cite{Gust-CMB-2011,Mortad-CMB-2011,mortad-CAOT-sum-normal} and
further bibliography cited therein.

In this section we want to start out ''from scratch'', assuming only
that $A$ and $B$ are densely defined closeable operators.  In other
words, we want to go beyond any assumptions of $B(H)$, positivity,
closed ranges, for any of our operators, inasmuch as that may be
possible.  Of course, we will need to make other, lesser,
assumptions, as we go along.

Early on, following results of J. von Neumann \cite{neumann}, M. A.
Naimark \cite{naim} (see also \cite{CHERNOFF} and \cite{SCHMUDG})
showed that even a densely defined closed symmetric operator $T$ may
have domain $D(T^2)$ not dense, even $D(T^2) = \{ 0\}$. Therefore we
shall normally assume that $D(AB)$ is dense, so that we may speak of
$(AB)^*$. That $D(AB)$ be dense is often apparent in applications,
e.g. where $D(AB)$ contains the smooth $C_0^\infty$ functions.

Let $A$ and $B$ be densely defined closeable operators such that $AB$ is
also densely defined.  Then $(AB)^* \supset B^* A^*$.  Then the domain
$D((AB)^*)$ is also dense iff $AB$ is closeable, and given that, then
$((AB)^*)^*$ exists and equals $\overline{AB}$ which has dense domain
$D(\overline{AB})\supset D(AB)$.  Thus $(AB)^{**}$ is a closed densely
defined operator and if $D(B^* A^*)$ is also dense we arrive at the two
inclusion relations
$$(B^* A^*)^* \supset (AB)^{**} = \overline{AB} \supset AB$$
and
$$(B^* A^*)^* \supset A^{**} B^{**} = \overline{A}~\overline{B} \supset
AB.$$
We remark that the distinction between $\overline{AB}$ and $\overline{A}
~\overline{B}$ is important and nontrivial so we want to keep it in mind
even when we are thinking of closed $A$ and $B$.

As all of the operators mentioned are densely defined, we may again
adjoint the above relations, arriving at, respectively,
$$(AB)^* = (\overline{AB})^* = (AB)^{***} \supset (B^* A^*)^{**} \supset
B^* A^* = (\overline{B})^* (\overline{A})^*$$
and
$$(AB)^* \supset (\overline{A}~\overline{B})^* = (A^{**} B^{**})^*
\supset B^* A^* = (\overline{B})^* (\overline{A})^*.$$
Doing this again serves as both check and possible refinement of these
considerations, from which we have, respectively,
$$\overline{AB} = (\overline{AB})^{**} \subset (B^* A^*)^{***} = (B^*
A^*)^* = (\overline{B^* A^*})^* = (\overline{\overline{B}^*
~\overline{A}^*})^*,$$
and
$$\overline{AB} \subset (\overline{A} ~\overline{B})^{**} =
(\overline{\overline{A} ~\overline{B}})^{**} = \overline{\overline{A}
~\overline{B}} \subset (B^* A^*)^*.$$

Let us summarize the relations obtained to this point.

\begin{lem}  Given densely defined  closeable $A$ and $B$ such
that $D(AB)$ is dense, we have $(AB)^* \supset B^* A^*$.  If also
$D(B^* A^*)$ is dense, we have all of the above relations.  In
particular, $\overline{AB}$ and $\overline{A} ~\overline{B}$ are
both included between $AB$ and $(B^* A^*)^*$.
\end{lem}

Of course one could see that $\overline{AB}$ and $\overline{A}
~\overline{B}$ are both sandwiched between $AB$ and $(B^* A^*)^*$
directly from
$$(B^* A^*)^* \supset A^{**} B^{**} = \overline{A} ~\overline{B},$$
i.e., from the general adjoint property, and from the general property
that the operator closure is always minimal, hence
$$(B^* A^*)^* \supset \overline{AB}.$$
But we wanted also the finer details of the Lemma available to us.

From the above it is clear that there are intimate relationships between
product adjoints and product closures and closures of products.  In
particular, we may ask the consequences of when equality may hold in the
above relations.  For example, suppose equality holds in the fundamental
adjoint operation
$$(B^* A^*)^* = A^{**} B^{**}.$$
If $A$ and $B$ be closed, that means $AB$ is closed.  Moreover, then
$(AB)^* = \overline{B^* A^*}$.  Should $B^* A^*$ also be closed, then the
fundamental adjoint relation $(AB)^* = B^* A^*$ obtains.  Similarly,
suppose equality holds in the other above obtained general adjoint
relation $(B^* A^*)^* \supset \overline{AB}$, i.e., that
$$(B^* A^*)^* = \overline{AB}.$$
Then since $(AB)^* = (\overline{AB})^*$ we have $(AB)^* = \overline{B^*
A^*}$, whether or not $A$ and $B$ are closed.  Should $B^* A^*$ also be
closed, then again the fundamental adjoint relation $(AB)^* = B^* A^*$
holds.  Let us summarize a version of these observations.

\begin{thm} Let $A$ and $B$ be densely defined closeable operators such
that $D(AB)$ and $D(B^* A^*)$ are both dense.  Then the product
adjoint equality $(B^* A^*)^* = A^{**} B^{**}$ necessitates that the
product $\overline{A} ~\overline{B}$ be closed and that
$(\overline{A} ~\overline{B})^* = (\overline{B^* A^*})$.  If
moreover $A$ and $B$ are closed, then $AB$ is necessarily closed.
If, further, the product $B^* A^*$ is closed, then the fundamental
adjoint relation
$$(AB)^* = B^* A^*$$
obtains.  Similarly, equality in the other adjoint relation above, $(B^*
A^*)^* = \overline{AB}$, along with closure of the product $B^* A^*$,
also yields the sought adjoint relation $(AB)^* = B^* A^*$.  Conversely,
necessary conditions for that fundamental adjoint relation include:  $B^*
A^*$ closed, and $\overline{AB} \supset \overline{A} ~\overline{B}$.
\end{thm}

\begin{rema}
Assuming $B^*A^*$ closed may not be merely dropped. For if $A$ and
$B$ are self-adjoint, and $B$ is the bounded inverse of $A$, then
$AB=I$ and $BA\subset I$. Hence $(AB)^*\neq BA$. An explicit example
would be:
\begin{exa}\label{main example}
Let $A$ be the unbounded operator defined by $Af(x)=(1+x^2)f(x)$ for
all $f\in L^2(\R)$ on the domain $D(A)=\{f\in L^2(\R):~x^2f\in
L^2(\R)\}$. Then it is clear that $A$ is self-adjoint and positive
on $D(A)$. Let $B$ be its inverse, i.e.
\[Bf(x)=\frac{1}{1+x^2}f(x),~\forall f\in L^2(\R).\]
\end{exa}
\end{rema}

We may recover the following well-known result (due to von-Neumann)
\begin{cor}
If $A$ is a densely defined closed operator, then $AA^*$ and $A^*A$
are self-adjoint.
\end{cor}

From the above analysis, one is naturally led to wonder:  could $B^*
A^*$ closed and densely defined provide a sufficient condition for
the sought $(AB)^* = B^* A^*$ fundamental adjoint product relation?
After all, if $A$ and $B$ are also closed, then if $(B^* A^*)^* =
A^{**} B^{**} = AB$, and adjointing again yields the sought $(AB)^*
= B^* A^*$.  The answer is:  yes, usually.  This follows from the
following fundamental result obtained by one of us long ago
\cite{Gust-1969-BAMS}.

\begin{thm} Let $T$ and $S$ be densely defined linear operators in a
Hilbert space.  If $S$ is closed and range $R(S)$ has finite
codimension, then $(TS)^* = S^* T^*$.
\end{thm}

\begin{cor} Let $A$ and $B$ be closed densely defined operators with
$D(AB)$ and $D(B^* A^*)$ dense.  Suppose also that $B^* A^*$ is
closed and densely defined and $R(A^*)$ is closed.  Then
$$(AB)^* = B^* A^*$$
unless it happens that the null space $N(A)$ is infinite
dimensional.
\end{cor}

The proof of Corollary 2 uses Theorem 2 with $T = B^*$ and $S = A^*$
and to note that
$$\text{codim} ~R(A^*) = \dim N(A)$$
by the Fundamental Theorem of Linear Algebra.  Then by Theorem 2,
$(B^* A^*)^* = AB$ and the result follows.

In Theorem 2 and Corollary 2, respectively, certainly $R(S)$ and
$R(A^*)$ have finite codimension, namely, zero codimension, if they
are the whole space.  By Banach's great theorem, that occurs iff
$S^*$ and $A$, respectively, have bounded inverses in $B (H)$.
Probably the following lemma has been noted elsewhere but we state
it here anyway because we will  want to use it later.

\begin{lem}\label{(AB)*=B*A* B INVERTIBLE!!} If $A$ and $B$ are densely defined and $B$ is closed with
inverse $B^{-1}$ in $B(H)$, then $(AB)^* = B^* A^*$.  In particular,
this is the case when $B$ is unitary.
\end{lem}

For more on the standard Fredholm theory, which we are going beyond
in this paper, see Remark 4.2.

\section{Unbounded Self-adjoint and Normal Products}
The important and immediately confirmable result that two
self-adjoint operators $A$ and $B$ in $B(H)$ have a self-adjoint
operator product $(AB)^* = AB$ iff they commute, $AB = BA$, will now
be examined for unbounded self-adjoint operators $A$ and $B$.  Here,
as in the preceding section, we wish to start ``from scratch'', to
seek more generality.

To avoid possible confusion, we would like to make precise certain
definitions here.

\begin{defn}  An unbounded operator $A$ is said to commute with an
unbounded operator $B$ if $AB \subset BA$.  $B$ commutes with $A$ if
$BA \subset AB$.  $A$ and $B$ are said to commute if $A$ commutes
with $B$ and $B$ commutes with $A$.  Then $AB = BA$.
\end{defn}

\begin{thm} Let $A = A^*$ and $B = B^*$ be two unbounded self-adjoint
operators such that $AB$ is densely defined.  If $A$ commutes with
$B$, then $AB$ is a symmetric operator and we have the relations
$$(AB)^* \supset \overline{BA} \supset \overline{AB} \supset AB.$$
Further, then $AB$ is self-adjoint if and only if
$$(AB)^* = \overline{BA} = \overline{AB} \equiv AB.$$
\end{thm}

\begin{proof}  All follows directly from the adjoint relation
$$(AB)^* \supset B^* A^* = \overline{B} ~\overline{A} = BA$$
and the assumed commutativity of $A$ with $B$.  We may additionally note
that adjointing the commutativity first gives another inclusion
$$(AB)^* \supset (BA)^* \supset \overline{AB} \supset AB$$
which features $\overline{AB}$ rather than the closures $\overline{B}$
and $\overline{A}$.
\end{proof}

What can we say in the reverse direction?

\begin{thm}  Let $A = A^*$ and $B = B^*$ be two unbounded self-adjoint
operators such that their product $AB$ is densely defined and
self-adjoint.  Then $B$ commutes with $A$.  If also $D(BA)$ is dense
and $BA$ is self-adjoint, then $A$ commutes with $B$.
\end{thm}

\begin{proof}  From the general relationship $(AB)^* \supset B^* A^*$ we
have from the assumption that $AB$ be self-adjoint that
$$AB = (AB)^* \supset B^* A^* = BA.$$
So necessarily $B$ commutes with $A$.  If $D(BA)$ is dense, then by
Theorem 3 above, $BA$ is symmetric.  Moreover,
$$(BA)^* = (B^* A^*)^* \supset A^{**} B^{**} = \overline{A} ~\overline{B}
= A^* B^* = AB.$$ Hence $A$ commutes with $B$ when $BA$ is
self-adjoint.
\end{proof}

Thus, if two unbounded self-adjoint operators $A$ and $B$ are such
that both operator products $AB$ and $BA$ are self-adjoint, then
they are in fact equal, $AB = BA$, and $A$ and $B$ commute just as
in the classical $A$ and $B$ bounded operators situation.  Moreover,
the adjoint product rule obtains:  $(AB)^* = B^* A^* = BA$.

In the converse direction, if two unbounded self-adjoint operators
$A$ and $B$ have both product domains $D(AB)$ and $D(BA)$ dense and
both commute in the unbounded sense, i.e., $AB \subset BA$ and $BA
\subset AB$, then of course they are equal, $AB = BA$, and thus they
commute in the classical strong sense.  But we may only conclude
that the product is a symmetric operator.  As to its closure, we may
note that from relations we used above, one has
$$AB \subset \overline{AB} \subset (AB)^* = (BA)^* \supset \overline{A}
~\overline{B} = AB,$$
reminiscent of the delicate relationship between the closure of the
product and the product of the closures that we identified in Section 2.

The considerations of Theorems 3 and 4 extend to arbitrary densely
defined closed operators $A$ if we restrict $B$ to be $A^*$.  Then
we have the established $AA^*$ theory, in which both $AA^*$ and $A^*
A$ are densely defined self-adjoint operators on dense domains
$D(AA^*)$ and $D(A^* A)$ which are domain cores for $A^*$ and $A$,
respectively.  The $AA^*$ theory is also a special case of Theorem
2.  The assumptions of dense domains are met.  The key assumption
$(B^* A^*)^* = A^{**} B^{**}$ becomes the self-adjointness $(AA^*)^*
= AA^*$.  The conclusion $(AB)^* = B^* A^*$ is exactly the same.
There is no need to have a finite dimensional nullspace $N(A)$ as in
Corollary 2.

Finally, we turn to normal operators $A$ (recall that an unbounded
operator $A$ is normal if it is closed and $AA^*=A^*A$). A natural
question within the context of this paper is: given $A$ and $B$
unbounded normal operators, what can one say about their product
$AB$?  This question has recently been studied, assuming some
boundedness, in a number of papers.  See for example
\cite{Mortad-IEOT-2009,Mortad-Demm-math,Mortad-product-unbounded-kaplansky}
and literature citations therein.  We state some of those results in
the following theorem.

Recall that if $A$ and $B$ are both bounded and commuting normal
operators, then one directly verifies, via the Fuglede theorem, that
their product $AB$ is normal. Such is not generally the case for
unbounded operators. A useful counterexample is Example 1. Indeed,
$B$ is bounded and positive (hence normal!). Observe that $BA\subset
AB=I$, strictly. So we see that even though we have strong
hypotheses, $BA$ is not guaranteed to remain normal! In fact, $BA$
is not even guaranteed to be closed.

However, we have

\begin{thm}[\cite{Mortad-Demm-math,Mortad-product-unbounded-kaplansky}]\label{BA normal A uniatry many results}  Let $A$ and $B$ be normal operators.
\begin{enumerate}
  \item Then $AB$ is normal under each of the following sufficient
conditions:
\begin{enumerate}
  \item $A$ is unitary and $A$ commutes with $B$, i.e., $AB \subset
BA$.
  \item $B$ is unitary and $B$ commutes with $A$, i.e., $BA \subset
AB$.
\end{enumerate}
  \item When $A$ is unitary, then $AB$ is normal iff $BA$ is normal.
\end{enumerate}
\end{thm}

We note here the following very simple proof of Part (2) of the
theorem: since $A$ is unitary, we may write
\[AB=A(BA)A^* \text{ and } BA=A^*(AB)A.\]
An unbounded operator unitarily equivalent to a normal one is normal
too. We actually do not need $B$ being normal, any $B$ will do. Of
course, an assumption of unitarity is exceedingly strong. Below we
will weaken that assumption.

Now, going back to Example 1, we see that $BA$ has a normal closure.
This has intrigued us to ask the following question: Under what
general circumstances does $BA$ have a normal closure? The answer is
given by the following result

\begin{thm}\label{normal closure}
Let $B$ a bounded normal operator and let $A$ be an unbounded normal
operator. If $B$ commutes with $A$, then $\overline{BA}$ is normal.
\end{thm}

We can prove the foregoing theorem using the spectral theorem for
normal (not necessarily bounded) operators, but we give a different
proof which relies on the following nice result (unfortunately, not
known to many!) by Devinatz-Nussbaum-von-Neumann:

\begin{thm}\label{Devintaz-Nussbaum-Neumann}[\cite{DevNussbaum-von-Neumann}]
Assume that $T_1$, $T_2$ and $T$ are self-adjoint operators such
that $T\subseteq T_1T_2$. Then $T=T_1T_2$.
\end{thm}

Now we prove Theorem \ref{normal closure}.

\begin{proof}
First we write all the hypotheses: $B$ commutes with $A$ means that
$BA\subset AB$. Then $B^*A^*\subset A^*B^*$ and by the Fuglede
theorem, $BA^*\subset A^*B$. Second, it is clear that $BA$ is
densely defined.

Now saying $\overline{BA}$ is normal amounts to saying that
$\overline{BA}^*$ is normal or just $(BA)^*$ is normal. Let us then
show that $(BA)^*=A^*B^*$ is normal. First it is clearly closed
thanks to the closedness of $A^*$ and the boundedness of $B^*$.
Next, we have on the one hand:

\[BB^*AA^*\subset BAB^*A^*\subset BAA^*B^*\subset (A^*B^*)^*A^*B^*.\]
Hence \[(A^*B^*)^*A^*B^*\subset AA^*BB^*=A^*AB^*B\] (where we have
taken adjoints and used the fact that $(A^*B^*)^*A^*B^*$ is
self-adjoint).

On the other hand:

\[A^*B^*(A^*B^*)^*\subset A^*B^*BA\subset A^*B^*AB\subset A^*AB^*B.\]

Setting $S=A^*B^*$, we see immediately that
\[S^*S\subset A^*AB^*B \text{ and } SS^*\subset A^*AB^*B.\]
Since $S^*S$, $A^*A$ and $B^*B$ are all self-adjoint, Theorem
\ref{Devintaz-Nussbaum-Neumann} applies and gives us
\[S^*S= A^*AB^*B \text{ and } SS^*= A^*AB^*B\]
so that
\[S^*S=SS^*,\]
proving the normality of $S$.
\end{proof}

Of course we have the following consequence which also improves Part
1 (a) of Theorem \ref{BA normal A uniatry many results}.

\begin{cor}
Let $B$ be a bounded and invertible normal operator and let $A$ be
an unbounded normal operator. If $B$ commutes with $A$, then $BA$ is
normal.
\end{cor}

\begin{proof}
Since $B$ is invertible, $BA$ is closed and so $BA=\overline{BA}$.
\end{proof}

Here is also an improvement of Part 1 (b) of Theorem \ref{BA normal
A uniatry many results}.

\begin{pro}
Let $A$ be a bounded and invertible normal operator and let $B$ be
an unbounded normal operator. If $A$ commutes with $B$, then $BA$ is
normal.
\end{pro}

\begin{proof}
First, $BA$ is closed by the boundedness of $A$ and the closedness
of $B$. Next, since $A$ is invertible, Lemma \ref{(AB)*=B*A* B
INVERTIBLE!!} yields $(BA)^*=A^*B^*$. Then we have by the Fuglede
theorem

\[(BA)^*BA=A^*B^*BA\subset B^* A^*BA\subset B^*BA^*A\]
and
\[BA(BA)^*=BAA^*B^*\subset BAB^*A^*\subset BB^*AA^*.\]
Since $B$ and $A$ are normal, Theorem
\ref{Devintaz-Nussbaum-Neumann} implies that $BA$ is normal.
\end{proof}

Very recently, the following theorem (among other interesting
results) was proved

\begin{thm}[\cite{Mortad-ex-Madani}]
Let $A$ and $B$ be two normal operators. Assume that $B$ is bounded.
If $BA=AB$, then $BA$ (and so $AB$) is normal.
\end{thm}

Now we ask: What can one say without assuming any boundedness on
either $A$ and $B$? Here we would like to present on a preliminary
basis some initial key findings. For simplicity, we assume that $A$
and $B$ are normal with dense ranges. Then (see Remark 4.3) an
operator is normal iff its polar factors commute. We will
concentrate on the needed behavior of these polar factors. Our goal
is to see necessary conditions for the normality of the unbounded
operator product $AB$.

Let $A = U|A| = |A|U$ be normal on $D(A) = D(|A|)$ and $B = V|B| =
|B|V$ be normal on $D(B) = D(|B|)$, where we have represented $A$
and $B$ in their polar factorizations and of course $U$ and $V$ are
unitary. Let $D(AB)$ be dense and $AB$ be closed.  For $AB$ to be
normal, $D((AB)^*)$ must equal $D(AB)$, and for simplicity in order
to fix the key ideas, we assume that here. Then we have
$$AB = U|A|V|B| = |A| UV|B| = U|A||B|V$$
among other possible expressions.  Using the last of these and Lemma
2 we have
$$(AB)^* = V^* (|A||B|)^* U^*.$$
Thus
$$(AB)^* (AB) = V^*[(|A||B|)^*(|A||B|)]V$$
and
$$(AB) (AB)^* = U[(|A||B|)((|A||B|)^*]U^*.$$

We note that if $(AB)^*AB$ is to equal $AB(AB)^*$, necessarily both
$(AB)^*AB$ and $AB(AB)^*$  must have $D((AB)^2)$ as their common
dense domain. From the above polar representations, then we must
have
$$V^* (T^*T)V = U(TT^*)U^*$$
where $T = |A||B|$ is the product of two nonnegative self-adjoint
operators.

Continuing, for $AB$ to be normal, we must have $||ABx||=
||(AB)^*x||$ on their common domain $D(AB)$. Squaring those leads to
the necessary condition that
\[<(T^*T-TT^*)U^*x,U^*x>=0\]
for the dense set $\{U^*x|~x\in D((AB)^2)\}$. When this dense set is
sufficient to conclude that in fact $T^*T=TT^*$, we have $T=|A||B|$
normal.

Because as a nonnegative self-adjoint operator $|B|$ satisfies the
spectral condition that $\sigma (|B|)) \cap \sigma (-|B|) = \{ 0\}$
of Theorem 6 of \cite{Mortad-PAMS2003}, we know that $|A||B|$ is
normal iff it is self-adjoint.

By our Theorem 4 above, that means that $|A|$ and $|B|$ must
commute.

We may summarize the conclusions of our discussions above as
follows.

\begin{thm}\label{main theorem gustafson normality}  For the product $AB$ of two unbounded normal
operators $A$ and $B$ with dense ranges to be normal, it is
necessary that $D(AB) = D((AB)^*)$ be dense, that $AB$ be closed,
that $D((AB)^2)$ be dense, that $(AB)^*AB=AB(AB)^*$ on common domain
$D((AB)^2)$, and that the polar factors product $T=|A||B|$ have
commutator vanish on the dense set $U^*(D((AB)^2))$. Thus the
essential necessary condition for $AB$ to be normal is that the
Hermitian polar factors of $A$ and $B$ commute.
\end{thm}

We remark that the conditions of Theorem \ref{main theorem gustafson
normality} are more than sufficient for $AB$ to be normal. For the
sufficiency, one needs only the assumptions stated prior to our
conclusion about the polar factors. Then the sufficiency proof is
essentially contained in Theorem 5.40 in \cite{Weidmann} by letting
$A_1=AB$ and $A_2=(AB)^*$ there.

We would like to make one further point about necessary conditions
on the polar factors here. To do so,  we recall the following
convenient terminology from scattering theory, e.g., see \cite{Amr}.

\vskip 1em {\bf Definition,} An operator $K$ is said to intertwine
two operators $N$ and $M$ if $KN \subset MK$.

The Fuglede-Putnam Theorem \cite{Con} or \cite{Putnam-book} (see
also \cite{MHM-Fug-Put-GMJ}, \cite{Mortad-Fuglede-Putnam-CAOT-2011}
and \cite{Stochel-asymmetric-Fuglede-2001} for more generalized
versions of this theorem) may be stated usefully in this terminology
as: If a bounded operator $K$ intertwines two unbounded normal
operators $N$ and $M$, then $K$ also intertwines $N^*$ and $M^*$.

Going back now to Theorem \ref{main theorem gustafson normality}, we
see that another necessary condition for $AB$ to be normal is the
intertwining
$$U^* V^* (T^* T) = (TT^*) U^* V^*$$
by the combined unitary $U^* V^*$ from the polar factorizations of
the combined Hermitian factors $T^* T$ and $TT^*$ if $AB$ is to be
normal.

This conclusion is seen to be consistent with the Fuglede-Putnam
Theorem because $T^* T$ and $TT^*$ are self-adjoint.

Thus we may state

\begin{cor}\label{corllary 3 after necssary theorem normal!!!!}  Given the necessary domain compatibilities, unbounded normal
operators $A$ and $B$ with dense ranges have normal operator product
only if their unitary factor adjoints product $U^* V^*$ intertwines
their Hermitian factor products as shown above.
\end{cor}

We comment that we stated Corollary \ref{corllary 3 after necssary
theorem normal!!!!} in terms of $U^* V^*$ intertwining $T^* T$ and
$TT^*$ as the more primitive, basic intertwining, even though we
know from Theorem 7 that $T^* T$ and $TT^*$ must necessarily be the
same, namely the (self-adjoint) operator $(|A||B|)^2$.  As such, the
intertwining becomes a full commutativity.

See Remark 4.4 for some further discussion of this problem.

\section{Dirac Operators, Additional Results, and Discussion}

\subsection{Remark}After we wrote the first draft of this paper,
independently F. Gesztesy forwarded to one of us (MHM) his recent
extensive work with J. A. Goldstein, H. Holden, G. Teschl
\cite{Gesztesy} treating abstract wave equations and Dirac type
operators. Key to their approach is consideration of operators for
which 0 belongs to the continuous spectrum $\sigma_c(A^*A)$. That is
also a context to which our results apply. Recall, eg see
\cite{Gust-1968-State-diagrams} and \cite{Gust-State-dia-1997} that
for a closed operator $T$, $0\in \sigma_c(T)$ iff $T$ is one-to-one
and $R(T)$ is properly dense.

Gesztesy et al \cite{Gesztesy} point out that in contrast to the
well-known long-standing important "folklore" interest in the
problem of when equality holds between $(TS)^*$ and $S^*T^*$,
situations which relate $\overline{ST}$ and
$\overline{S}~\overline{T}$ have been much less studied. We
certainly agree. The following result is proved in \cite{Gesztesy}.
We use our notation.

\begin{pro}(\cite{Gesztesy})\hfill
\begin{enumerate}
  \item If $\overline{B}$ is in $B(H)$ and $A$ is closed, then
  $\overline{AB}\subset A\overline{B}$.
  \item If $\overline{B^{-1}}\in B(H)$, $A$ on $D(A)\cap R(B)$ is closeable,
  and $A$ on $D(A)\cap R(\overline{B})$ is contained in the closure
  of $A$ on $D(A)\cap R(B)$, and  $AB$ is closeable, then $A\overline{B}\subset \overline{AB}$.
\end{enumerate}
\end{pro}

The importance of the distinction between $\overline{AB}$ and
$\overline{A}~\overline{B}$ was impressed upon one of us (KG) long
ago (in our unpublished work [[3]] discussed in our paper
\cite{Gust-CMB-2011}, and in\cite{RKNG-1983} and
\cite{Gustafson-Rejto} too) by the investigations of Schmincke
\cite{Schmincke} of self-adjoint extensions of Dirac operators. The
following result is stated without proof in \cite{Schmincke}.

\begin{pro}(\cite{Schmincke}) Let $A$, $B$ and $AB$ be closeable operators.
Then
\begin{enumerate}
  \item If $B$ is relatively bounded with respect to $AB$, then
  $\overline{AB}\subset \overline{A}~\overline{B}$.
  \item If $B^{-1}$ is closeable and relatively bounded with respect
  to $A$, then $\overline{A}~\overline{B}\subset \overline{AB}$.
\end{enumerate}
\end{pro}

As is evident, \cite{Schmincke} criteria are more general than
\cite{Gesztesy}, allowing $B$ and $B^{-1}$ to be unbounded. We would
like to look more closely at these propositions from the framework
of this paper.

We recall for the reader's convenience the following rather general
sufficient criteria for relative boundedness of one operator to
another, sometimes called H\"{o}rmander's Theorem. See, e.g.,
\cite{Gold-Book-1966} or \cite{Yosida}.
\begin{thm}\label{hormander theorem??}
If $T$ is closed and $B$ is closeable with $D(T)\subset D(B)$, then
$B$ is $T$-bounded.
\end{thm}

We start with the following lemma, which brings our approach and
results in Section 2 to bear on these key Propositions 2 and 3 above
(called Lemmas 2 in both \cite{Gesztesy} and \cite{Schmincke}) used
in those respective analyses of Dirac operators.

\begin{lem}
Let $A$, $B$ and $AB$ all be densely defined and closeable. Then
\begin{enumerate}
  \item If $D(B^*A^*)$ is dense, we have $\overline{AB}\subset
  \overline{\overline{A}~\overline{B}}$.
  \item In particular, if $\overline{B}$ is in $B(H)$ and $A$ is
  closed, then $\overline{AB}\subset A\overline{B}$.
  \item Without $D(B^*A^*)$ necessarily dense, we may state
  generally, in the opposite direction, that $\overline{AB}\subset
  \overline{A}~\overline{B}$ implies that $\overline{B}$ is
  relatively bounded with respect to $\overline{AB}$, and hence, $B$
  is relatively $AB$-bounded.
  \item Also, generally, $\overline{A}~\overline{B}\subset
  \overline{AB}$ implies that $\overline{A}~\overline{B}$ is
  closeable and
  $\overline{\overline{A}~\overline{B}}=\overline{AB}$.
\end{enumerate}
\end{lem}

\begin{proof}
We established (1) on proving Lemma 1 in Section 2. Then (2) follows
from (1) because $B^*$ is in $B(H)$ and hence $D(B^*A^*)=D(A^*)$ is
dense. Turning to (3), suppose that $\overline{AB}\subset
\overline{A}~\overline{B}$. Then we have the situation of Theorem 10
above, namely
\[D(\overline{AB})\subset D(\overline{A}~\overline{B})\subset D(\overline{B})\]
for $\overline{AB}$ closed and $\overline{B}$ closeable. Thus
$\overline{B}$ is $\overline{AB}$-bounded. Thus $B$ is $AB$-bounded,
since always $D(AB)\subset D(B)$. As to (4), suppose $AB\subset
\overline{A}~\overline{B}\subset \overline{AB}$. Generally, whenever
you have a sandwiching $S\subset T\subset \overline{S}$, then
$T^*=S^*$ and $T$ is closeable. Let us adjoint twice as in Section
2, giving us the relations verifying (4).
\[AB\subset \overline{A}~\overline{B}\subset \overline{AB}\]
\[(AB)^*=(\overline{A}~\overline{B})^*=(\overline{AB})^*=(AB)^*\]
\[(AB)^{**}=\overline{AB}=(\overline{A}~\overline{B})^{**}=\overline{\overline{A}~\overline{B}}.\]

\end{proof}
We note that (2) of Lemma 3 is Proposition 2 (1) which is thus seen
to be a special case of our results of Section 2. We note that (3)
of Lemma 3 reverses the direction of Proposition 3 (1) and shows
that the sufficient condition assumed there is actually necessary.
And (4) of Lemma 3 provides a necessary condition for both
Proposition 2 (2) and Proposition 3 (2).

We were not able to provide a proof of the sufficiency of
Proposition 3 (1) in the generality with which it is stated. A
formal argument would be the following. If $B$ is $AB$-bounded put
the $\overline{AB}$ operator norm on $D(\overline{AB})$. Then $B$,
being $AB$-bounded, is bounded on that newly-normed space and has a
bounded extension $\overline{B}$ in $B(\overline{AB}, H)$. Then
$\overline{A}~\overline{B}$ is in fact closed by the closedness of
$\overline{A}$ and boundedness of $\overline{B}$. Anytime
$\overline{A}~\overline{B}$ is closed, necessarily
$\overline{AB}\subset \overline{A}~\overline{B}$ by the minimality
of the closure property. But there is is an inconsistency. The
closure on the left is in the original Hilbert space whereas the one
on the right combines closures in two different norms.

We are able, however, to obtain a result a bit weaker than that
claimed in Proposition 3 (1), i.e. that claimed in \cite{Schmincke}
(Lemma 2 (i)).

\begin{lem}\label{lemma 4}Let $A$, $B$ and $AB$ be closeable. If $\overline{B}$ is
$\overline{A}~\overline{B}$-bounded, then
$\overline{A}~\overline{B}$ is closed and $\overline{AB}\subset
\overline{A}~\overline{B}$.
\end{lem}

\begin{proof}
Let $x_n\in D(\overline{A}~\overline{B})\subset D(\overline{B})$
converge to $x$ and $\overline{A}~\overline{B}x_n=y_n$ converge to
$y$. Then by the assumed relative boundedness,
\[\|\overline{B}(x_n-x_m)\|\leq a\|x_n-x_m\|+b\|\overline{A}~\overline{B}(x_n-x_m)\|,\]
and $(Bx_n)$ is a Cauchy sequence. Thus $\overline{B}x_n$ converges
to some $z$. $\overline{B}$ being closed guarantees $x\in
D(\overline{B})$ and $\overline{B}x=z$. The closedness of
$\overline{A}$ then puts $z\in D(\overline{A})$, which in turn puts
$x\in D(\overline{A}~\overline{B})$. Thus
$\overline{A}~\overline{B}x=\overline{A}z=y$ and
$\overline{A}~\overline{B}$ is closed.

In the other direction, we may also note that the conclusion of
Lemma \ref{lemma 4} requires necessarily that $\overline{B}$ be
$\overline{A}~\overline{B}$-bounded. This follows from
$\overline{A}~\overline{B}$ being closed and
$D(\overline{A}~\overline{B})\subset D(\overline{B})$ and Theorem
\ref{hormander theorem??}.
\end{proof}

We turn next to parts (2) of Propositions 2 and 3. The proof of
Proposition 2 (2) given in \cite{Gesztesy} seems clear. The
attention to domains detail is gratifying and is akin to our
considerations in \cite{Gust-CMB-2011} and \cite{Mortad-CMB-2011}.

On the other hand, no proof of Proposition 3 (2) was given in
(\cite{Schmincke}, Lemma 2 (ii)). Rather than jumping in to attempt
a proof here, we are more disposed to first look more closely at its
hypotheses to understand them better, and compare them in some
detail to the hypotheses of Proposition 2 (2).

Specifically, we are able to provide some clarification by use of
our state diagrams of \cite{Gust-1968-State-diagrams} and
\cite{Gust-State-dia-1997}. Such are very useful as a summary of
those possible configurations between an operator $T$ and its
adjoint $T^*$ as to their ranges and inverses. One lets I, II, III
denote respectively $R(T)$ the whole Hilbert space, $R(T)$ properly
dense, $R(T)$ closure not dense, and the same for $T^*$. One lets
1,2,3 denote $T^{-1}$ bounded, $T$ just $1-1$, or $T$ not $1-1$ and
the same for $T^*$. Then of the 81 states formed by matching the 9
for $T$ versus the 9 for $T^*$, only 16 can actually occur. See
\cite{Gust-1968-State-diagrams,Gust-State-dia-1997} and especially
(\cite{Gust-1968-State-diagrams}, Figure 1). For densely defined
operators in a Hilbert space, only 13 states can occur. As you
demand closed operators, or $1-1$ operators, the possible states
decrease. For example, a self-adjoint operator $T$ can only be in
states $I_1I_1$, $II_2II_2$ or $III_3III_3$. See
\cite{Gust-1968-State-diagrams,Gust-State-dia-1997} and the
references and history therein for more details.

We also will use some well-known facts about closeable operators.
First, for $B$ closeable and $1-1$, $B^{-1}$ is closeable if and
only if $\overline{B}$ is $1-1$, and then
$(\overline{B})^{-1}=\overline{B^{-1}}$. See (\cite{Weidmann}, p.
90). Second, for $B$ with dense range and $1-1$, $B^*$ is also $1-1$
and $(B^*)^{-1}=(B^{-1})^*$. See (\cite{Weidmann}, p. 71).

We are particularly interested in the hypotheses of Proposition 3
(2). By systematic process of elimination through the state
diagrams, one may ascertain that for $B$ closeable, $1-1$ with
closeable inverse $B^{-1}$, $\overline{B}$ can only occur in 4
states: $I_1I_1$, $II_2II_2$ , $III_1I_3$, $III_2II_3$. In
particular, $R(B^*)=R((\overline{B})^*)$ is necessarily dense. But
then assuming that $B^{-1}$ is $A$-bounded means $D(A)\subset
D(B^{-1})=R(B)$. So $R(B)$ is necessarily dense and so is
$R(\overline{B})$ and one sees that $B$ can only occur in 3 states
$I_1I_1, II_1I_1, II_2II_2$ and $B^*$ is now $1-1$. Similarly,
$\overline{B}$ can have only states $I_1I_1$ and $II_2II_2$. One
does not need $D(A)\subset R(B)$ for this, just $R(B)$ dense.

In contrast, the hypotheses of Proposition 2(2), in particular, that
$\overline{B^{-1}}\in B(H)$, means that $B$ is closeable, $1-1$,
with closeable inverse as in Proposition 3 (2), but also $B^{-1}$ is
bounded and $R(B)$ is dense. Thus $B$ can only have states $I_1I_1$,
$II_1I_1$ so $R(B^*)$ is the whole space and $B^*$ has a bounded
inverse. Moreover $\overline{B}$ can only be in state $I_1I_1$.

This analysis enables a better comparison of the hypotheses of
Propositions 2 and 3, parts (2). The main distinction seems to be
that Proposition 2 needed to assume a bounded everywhere defined
inverse $\overline{B}^{-1}$ and then work on $D(A)\cap R(B)$,
whereas Proposition 3 allowed $B^{-1}$ unbounded but needed
$D(A)\subset R(B)$. In both instances one has $R(B)$ and
$R(\overline{B})$ both dense but in the former case
$R(\overline{B})$ is the whole space.

We leave any direct rigorous proof of Proposition 3 (2) as stated to
others. However, as the following variation on that situation shows,
if one interprets $B^{-1}$ being $A$-bounded to mean
$\overline{B}^{-1}$ is bounded, then you see immediately that you
are closer to the hypotheses of part (2) of Proposition 2. Then
proceeding by formal arguments similar to what we did above for part
(1) of Proposition 3, versions of part (2) may be concluded. But
$\overline{B}^{-1}$ is in $B(\overline{A},H)$ and so different
closure topologies are involved.

Thus, let $A$, $B$ and $AB$ be closeable as in Section 2, but assume
also that $A$ and $B$ are both $1-1$ with closeable inverses and
dense ranges. This is often the situation when one is working over
the continuous spectrum, e.g., see
\cite{Gesztesy,Gust-State-dia-1997,Gust-CMB-2011,Schmincke}. Our
idea is to approach hypotheses involving inverses, as in parts (2)
of Propositions 2 and 3, by combining inversing, adjointing, and
closing. Remember that $T^*=\overline{T}^*$ for a closeable operator
$T$.

We first invert then adjoint:
\[AB\subset \overline{A}~\overline{B}, B^{-1}A^{-1}=(AB)^{-1}\supset (\overline{A}~\overline{B})^{-1}=\overline{B}^{-1}~\overline{A}^{-1},\]
\[(A^{-1})^{*}(B^{-1})^*\subset ((AB)^{-1})^*\subset (\overline{B}^{-1}\overline{A}^{-1})^*=(\overline{A}^{-1})^*(\overline{B}^{-1})^*,\]
where we went ahead and assumed $\overline{B}^{-1}$ bounded for the
last equality. This expression can be rewritten as
\[(A^*)^{-1}(B^*)^{-1}\subset ((AB)^*)^{-1}\subset (\overline{A}^*)^{-1}(\overline{B}^*)^{-1}\]
which we may invert to

\[B^*A^*\supset (AB)^*\supset (\overline{B})^*(\overline{A})^*.\]
Then adjointing this and using closures, we arrive at the inclusion
$\overline{A}~\overline{B}\subset \overline{AB}$ sought in parts (2)
of Propositions 2 and 3. Using $D(B^*A^*)$ dense and Lemma 3 (1), we
may conclude that $\overline{A}~\overline{B}$ is closeable.
Summarizing, we have
\[\overline{A}~\overline{B}\subset \overline{\overline{A}~\overline{B}}\subset \overline{AB}\subset (B^*A^*)^*.\]

The right-most term could seem to be of possible interest in some
situations. From the $B^{-1}$ bounded assumption and the state
diagrams we know that $R(B^*)$ is the whole space. When it happens
that $B$ also is bounded on a dense domain, then all four operators
are equal.

\subsection{Remark}
If one assumes that both $T$ and its range $R(T)$ are closed, then
one has the established normally solvable operator theory. If
additionally index$(T)=\dim N(T)-\dim H/R(T)$ is well-defined, one
has the semi-Fredholm operator theory, and if index$(T)$ is finite,
the Fredholm theory for both bounded and unbounded operators. See
\cite{Gold-Book-1966} and \cite{Kato-Book} among others, and the
cited papers \cite{Castr-Gold}, \cite{Gust-1969-BAMS},
\cite{Gust-1969}
 and \cite{Sch-1970}. See also the recent paper
\cite{Azz-Dje-Mess}. That is not the case when one is on the
continuous spectrum.

We would like to point out the important result in
\cite{Castr-Gold}, Theorem 2 there: For given closed $S$, if $(TS)^*
= S^* T^*$ for all closed $T$, then necessarily $R(S)$ must have
finite codimension. That means $R(S)$ is also closed. Therefore, we
found it more convenient in Corollary 1 to state the codimension
condition in terms of the null space $N(A)$.
\subsection{Remark}
Generally a closed operator $T$ has a unique polar decomposition
$T=W|T|$ where $W$ is a maximal partial isometry from
$\overline{R(|T|)}$ to $\overline{R(T)}$ and with null space
$N(W)=N(|T|)$. See  \cite{Conway-book-AMS-GSM-2000},
 \cite{Halmos}, \cite{Kato-Book} and \cite{Weidmann}. Then one asks:
When do the isometric and Hermitian factors commute? For finite
dimensional matrices, the answer is iff $T$ is a normal matrix.
However, for infinite dimensional $T$ in $B(H)$, $W|T|=|T|W$ iff $T$
is quasi-normal: $T$ commutes with $T^*T$ (see
\cite{Conway-book-AMS-GSM-2000} or \cite{Halmos}). An extension of
this fact to $T$ unbounded quasi-normal is given in
\cite{MAjdak-2007}.

We did not wish to go beyond issues of normality of operator
products in this paper. Therefore as stated before Theorem \ref{main
theorem gustafson normality} in Section 3, we assumed that $T$ has
dense range $R(T)$. This rules out the possibility of quasi-normal
operator products unless they are normal. See the state diagram
papers \cite{Gust-1968-State-diagrams} and
\cite{Gust-State-dia-1997}.

To be precise here, we state the following

\begin{pro}
Let $T$ be quasi-normal with dense range $R(T)$. Then $T$ is normal.
\end{pro}

\begin{proof}
Because in its polar decomposition $T=W|T|$, the isometric factor is
unitary, we have (using Lemma \ref{(AB)*=B*A* B INVERTIBLE!!}) that
$T^*T=|T|^2$, and using the commutativity of $W$ with $|T|$ due to
quasi-normality, that
\[TT^*=W|T|^2W^*=|T|^2\]
\end{proof}

Recall (see \cite{Gust-1968-State-diagrams} and
\cite{Gust-State-dia-1997}) that for a closed densely defined
operator $T$ in a Hilbert space, the dense range $R(T)$ condition is
equivalent to an assumption that $T^*$ is one-to-one.

\^{O}ta \cite{Ota-2002} explores a class of $q$-quasi-normal
operators as a generalization of the concept of quasi-normal
operator, and gets results also of interest for quasi-normal and
general operators and especially for contractions.

\subsection{Remark}
An initial investigation of the general problem of when unbounded
normal operators $A$ and $B$ have normal product $AB$ was made long
ago by Devinatz-Nussbaum \cite{DevNussbaum}. Following the
convention of Riesz-Sz.Nagy \cite{Riesz-Nagy}, $A$ and $B$ are said
to permute if all of their canonical resolutions of the identity
commute. That is of course a very strong commutativity assumption.
Digressing a bit, we say that in fact many authors call the
permutability a strong commutativity or just commutativity. For
further points, we refer the interested reader to \cite{RS1} where
this question is treated with great detail together with the famous
Nelson example (see \cite{Nelson}). The whole section was then
called "Formal manipulation is a touchy business: Nelson's example".

Going back to \cite{DevNussbaum}, Property $P$ is then defined as
follows: $A$, $B$ and $AB$ are said to have property $P$ if they are
all normal and $AB=BA$. In other words, the normality of the
products $AB$ and $BA$ is assumed and also their full commutativity.
The main result of Devinatz-Nussbaum \cite{DevNussbaum} may then be
stated as follows, in our notation.

\begin{thm}\label{Devintaz-Nussbaum}
If $A$, $B$ and $AB$ have property $P$, then $A$ and $B$ permute.
\end{thm}

In other words, a necessary condition for $AB$ and $BA$ to both be
normal and furthermore to commute is that their detailed spectral
families must also all commute.

In our investigation in Theorem \ref{main theorem gustafson
normality}, we a priori wanted to start with no commutativity
assumptions. But we were led to an essential necessary condition
that their Hermitian polar factors commute.

On the other hand, if one is willing to assume sufficient
commutativity, one may obtain a number of sufficiency conditions.
For example:

\begin{pro}(\cite{mortad-sum-newest})\label{AB normal A B Invertible}
Let $A$ and $B$ be two unbounded invertible normal operators. If
$BA=AB$, then $BA$ (and $AB$) is normal.
\end{pro}

The proof uses the following version of the Fuglede theorem.
\begin{thm}(\cite{mortad-sum-newest})\label{Fuglede-NEW-ALL-UNBD}
Let $A$ and $B$ be two unbounded normal and invertible operators.
Then
\[AB=BA\Longrightarrow AB^*=B^*A \text{ and } BA^*=A^*B.\]
\end{thm}

Next, as an important consequence of Theorem \ref{Devintaz-Nussbaum}
we have:

\begin{cor}
Let $A$ and $B$ be two unbounded invertible normal operators. If
$BA=AB$, then $A$ and $B$ permute.
\end{cor}

With Theorem \ref{Devintaz-Nussbaum-Neumann}, we may prove a
stronger result. We still need a version of Theorem
\ref{Fuglede-NEW-ALL-UNBD}. We have:

\begin{thm}\label{Fuglede-NEW-ALL-UNBD-Only one is invertible}
Let $A$ and $B$ be two unbounded normal operators. If $B$ is
invertible, then
\[AB=BA\Longrightarrow A^*B\subset BA^* \text{ and hence } AB^*\subset B^*A.\]
\end{thm}

We include a proof for the reader's convenience.

\begin{proof}
Since $B$ is invertible, we may write
\[AB=BA \Longrightarrow B^{-1}A\subset AB^{-1}.\]
Since $B^{-1}$ is bounded and $A$ is unbounded and normal, by the
classic Fuglede theorem we have
\[B^{-1}A\subset AB^{-1}\Longrightarrow B^{-1}A^*\subset A^*B^{-1}.\]
Therefore,
\[A^*B\subset BA^*.\]
But $B$ is invertible, so by Lemma \ref{(AB)*=B*A* B INVERTIBLE!!}
we may obtain
\[AB^*\subset(BA^*)^*\subset (A^*B)^*=B^*A.\]
\end{proof}

Accordingly, we may state and prove the following

\begin{thm}\label{last result AB=BA normal B invertible!!!!!}
Let $A$ and $B$ be two unbounded normal operators such that only $B$
is invertible. If $AB=BA$, then $AB$ is normal.
\end{thm}

\begin{proof}
The proof is based on Lemma \ref{(AB)*=B*A* B INVERTIBLE!!} and
Theorem \ref{Fuglede-NEW-ALL-UNBD-Only one is invertible}. The use
of each of these results will be clear to the reader. The closedness
of $AB$ follows from that of $BA$ for $B$ is invertible and $A$ is
closed.

We have
\[AB(AB)^*=ABB^*A^*=AB^*BA^*\subset B^*ABA^*=B^*BAA^*.\]
By Theorem \ref{Devintaz-Nussbaum-Neumann}, we obtain

\[AB(AB)^*=B^*BAA^*.\]
We also have

\[(AB)^*AB=B^*A^*AB=B^*A^*BA\subset B^*BA^*A.\]
Applying Theorem \ref{Devintaz-Nussbaum-Neumann} again implies that
\[(AB)^*AB=B^*BA^*A.\]
Therefore, $AB$ is normal.
\end{proof}

Combining Theorem \ref{last result AB=BA normal B invertible!!!!!}
with Theorem \ref{Devintaz-Nussbaum} yields the following
consequence.

\begin{cor}
Let $A$ and $B$ be two unbounded normal operators such that $B$ is
invertible. If $BA=AB$, then $A$ and $B$ permute.
\end{cor}

\subsection{Remark}
One encounters many interesting technicalities when working with
unbounded operators. In this paper we looked at certain intricacies
in the interactions of the two operations of adjointing and closing
unbounded operator products. However, for the readers convenience
and future reference, we may present the closure situation, seen
separately, in a rather simplified way.

Suppose $A$, $B$ and $AB$ are all densely defined and closeable in a
Hilbert space. Suppose moreover that $\overline{A}~\overline{B}$ is
closeable. Then $AB\subset\overline{A}~\overline{B}\subset
\overline{\overline{A}~\overline{B}}$. Since $\overline{AB}$ is the
minimal closure operation, one has $\overline{AB}\subset
\overline{\overline{A}~\overline{B}}$. That is the generic situation
to keep in mind.

This key requirement that $\overline{A}~\overline{B}$ is closeable
is equivalent to $D[(\overline{A}~\overline{B})^*]$ being dense.
Since $(\overline{A}~\overline{B})^*\supset B^*A^*$, in Section 2 we
assumed for simplicity that $D(B^*A^*)$ be dense, which might be
easier to establish in applications.

In the case that $\overline{A}~\overline{B}$ is not closeable,
clearly one cannot have $\overline{A}~\overline{B}\subset
\overline{AB}$.

Finally, all of the considerations of this paper apply to
$\overline{A}~\overline{B}$ replaced with either $A\overline{B}$ or
$\overline{A}B$. We omit the details.

\section{Conclusions}

With an eye toward questions about products of unbounded
self-adjoint operators and products of unbounded normal operators
and certain operator-product issues for Dirac operators, we
established general results for operator closures
$\overline{A}~\overline{B}$ and $\overline{AB}$ for unbounded
closeable operators $A$ and $B$ in a Hilbert space. These results
are naturally tied to questions about operator product adjoints
$(AB)^*$. No commutativity of operators was assumed.

Then for conditions for operator products $AB$ to be self-adjoint or
normal, commutativity is unavoidable. We established both sufficient
and necessary conditions. In particular, for a product of two
unbounded normal operators $AB$ to be normal, it is necessary that
their polar factors $|A|$ and $|B|$ commute.

We also clarified certain fundamental operator-product issues for
Dirac operators.

\section*{Acknowledgement}
We very much appreciate Fritz Gesztesy forwarding his joint work
\cite{Gesztesy} to us during the course of our investigations.

\end{document}